\documentclass[12pt]{article}

\input amssym.def
\input amssym

\usepackage{theorem}

\newcommand{\CC}{{\Bbb C}}

\newcommand{\ZZ}{{\Bbb Z}}

\newcommand{\fg}{{\frak g}}
\newcommand{\fs}{{\frak s}}

\newcommand{\req}[1]{(\ref{eq:#1})}

\newcommand{\Aut}{\mathop{\rm Aut}}
\newcommand{\pdiff}{\mathop{\rm PDiff}}
\newcommand{\Res}{\mathop{\rm Res}}

\newcommand{\resp}{\hbox{\rm resp.}\ }

\newcommand{\subsym}[2]{%
  \mathrel{\raise-2ex\hbox{$\stackrel {\hbox{\normalsize $#1$}}%
      {\hbox{\scriptsize $#2$}}$}}}

\makeatletter
\@addtoreset{equation}{section}
\makeatother

\newcounter{saveeq}
\newenvironment{mathletters}%
{ \refstepcounter{equation}
  \setcounter{saveeq}{\value{equation}}
  \setcounter{equation}{0}
  }%
{ \setcounter{equation}{\value{saveeq}}%
  }

\newtheorem{proposition}{Proposition}[section]
\newtheorem{lemma}{Lemma}[section]
\newtheorem{definition}{Definition}[section]
\newtheorem{theorem}{Theorem}[section]
\newtheorem{corollary}{Corollary}[section]

{\theorembodyfont{\rmfamily}
  \newtheorem{example}{Example}[section]
  \newtheorem{remark}{Remark}[section]
}
\newenvironment{proof}[1][Proof]
           {\medbreak\noindent {\em #1: \enspace}}
           {\hfill \openbox \par \medbreak}

\newcommand{\openbox}{\leavevmode
  \hbox to.77778em{%
  \hfil\vrule
  \vbox to.675em{\hrule width.6em\vfil\hrule}%
  \vrule\hfil}}

\newcommand{\alphalist}{% changes enumerate 1st level to a)...z)
  \renewcommand{\theenumi}{\rm\alph{enumi}}%
  \renewcommand{\labelenumi}{\theenumi)}%
}

\newcommand{\alphaparenlist}{% changes enumerate 1st level to (a)...(z)
  \renewcommand{\theenumi}{\alph{enumi}}%
  \renewcommand{\labelenumi}{(\theenumi)}%
}

\begin{document}

\newcounter{bean}
\newenvironment{deflist}[1]%
    {
      \begin{list}{\bf #1\arabic{bean}}
         {\usecounter{bean}
              \setcounter{bean}{-1}
          \labelsep=1em
          \settowidth{\labelwidth}{#1\thebean:}
          \addtolength{\labelwidth}{1.1ex}
          \leftmargin=\labelwidth
          \addtolength{\leftmargin}{\labelsep} }
          \renewcommand{\makelabel}[1]{\upshape(##1)\hfil}
    }
    {\end{list}}

\newenvironment{optenum}[1]%
    {
    \begin{list}{}
        { \labelsep=1em
          \settowidth{\labelwidth}{#1}
          \addtolength{\labelwidth}{1.1ex}
          \leftmargin=\labelwidth
          \addtolength{\leftmargin}{\labelsep}
          \renewcommand{\makelabel}[1]{##1\hfil}
          }}
    {\end{list}}

\textwidth=16.3cm
\hoffset-1.5cm
\title{
\rightline{\small ITEP-TH-28/97}
\bigskip
\bigskip
\bigskip
$\Gamma$-Conformal Algebras}
\author{Maria I. Golenishcheva-Kutuzova\thanks{%
Institut Girard Desarques, URA 746, Universit\'{e} Lyon-1,%% \hfill\break
      $<$kutuzova@geometrie.univ-lyon1.fr$>$ %%  \hfill\break
      and %% \hfill\break
      International Institute of Nonlinear Science (Moscow)
    } ,
\and
    %%% second author
    Victor G. Kac\thanks{%
    Department of Mathematics, MIT,
    Cambridge, MA 02139, USA, %% \hfill\break
    $<$kac@math.mit.edu$>$
    }}
\date{}

\maketitle

\begin{abstract}
$\Gamma$-conformal algebra is an axiomatic description of the operator product
expansion of chiral fields with simple poles at finitely many points.
We classify these algebras and their representations in terms of Lie
algebras and their representations with an action of the group $\Gamma$. To
every $\Gamma$-conformal algebra and a character of $\Gamma$ we associate
a Lie algebra  generated by fields with the OPE with simple poles.
Examples include twisted
affine Kac-Moody algebras, the sin algebra
(which is a ``$\Gamma$-conformal'' analodue of the general linear
algebra) and its analogues, the algebra of pseudodifferential operators
on the circle, etc.

\end{abstract}

\setcounter{section}{-1}
\section{Introduction}
\label{sec:0}

In the past years there has been a number of papers, where the
effect of splitting of a pole of order $N$ in the operator
product expansion (OPE) occurs.  Examples include
representations of twisted affine Kac-Moody algebras
\cite{K-K-L-W}, of the Lie algebras of Quantum torus \cite{G-K-L},
 of Quantum affine algebras \cite{F-J}, and of
the central extension of the double of Yangian \cite{Kh-L-P}.
These constructions found many applications in the theory of
integrable systems and exactly solvable models in quantum field
theory.

The simplest motivation to consider the OPE with simple poles,
but in the shifted points is the following.  Consider the OPE of
two local fields $a (z)$ and $b (z)$:
\begin{equation}
  \label{eq:0.1}
  a (z) b(w) \sim \sum^{N-1}_{j=0} \frac{c^j (w)}{(z - w)^{j + 1}}
\end{equation}
and replace $(z - w)^n$ by its $q$-analogue $(z - w)^n_q = (z -
w) (z - q w) \cdots (z - q^{n - 1} w)$.  The result of this
$q$-deformation is an OPE with only simple poles.  A natural
generalization is to consider $\Gamma$-deformation, where
$\Gamma$ is a subgroup of the multiplication group of non-zero
complex numbers $\CC^\times$.  Namely, for an $n$-element subset $S$ of
$\Gamma$ we let
\begin{displaymath}
  (z - w)^n_S = \prod_{\alpha \in S} (z - \alpha w)
\end{displaymath}
and consider a collection of fields satisfying the condition of
$\Gamma$-locality, i.e.
\begin{equation}
  \label{eq:0.2}
  (z - w)^n_S \left[ a (z), b(w) \right] = 0
\end{equation}

The study of algebraic properties of $\Gamma$-locality leads us to
a formal definition of a $\Gamma$-conformal algebra given below.
Property \req{0.2} is equivalent to the following form of the OPE
(Proposition~1.1):
\begin{equation}
  \label{eq:0.3}
  a (z) b (w) \sim \sum_{\alpha \in S}
  \frac{c^\alpha (w)}{ z - \alpha w } \, .
\end{equation}
We view the field $c^\alpha (w)$ as an ``$\alpha$-product'' of
the fields $a (z)$ and $b (z)$:
\begin{equation}
  \label{eq:0.4}
  c^\alpha (w) = (a_{(\alpha)} b) (w) \, .
\end{equation}
Introduce operators $T_\alpha (\alpha \in \CC^\times)$ on the space of
fields by
\begin{equation}
  \label{eq:0.5}
  T_\alpha a (z) = \alpha a (\alpha z) \, .
\end{equation}
Let $R$ be a space over $\CC$ of $\Gamma$-local fields which is
closed under products \req{0.4} and invariant under operators
\req{0.5} for all $\alpha \in \Gamma$  Then one can show
(Proposition~2.1) that the following properties hold, where $a,
b, c \in R$,
\hbox{$\alpha, \beta \in \Gamma$}:

\begin{deflist}{C}

\item
  \label{item:C0}
  $a_{(\alpha)} b = 0$ for all but finitely many $\alpha \in \Gamma$,

\item
  \label{item:C1}
  $(T_\alpha a)_{(\beta)} b = a_{(\beta \alpha)} b$,

\item
  \label{item:C2}
  $a_{(\alpha)} b = - T_\alpha (b_{(\alpha^{-1})} a)$,

\item
  \label{item:C3}
  $[a_{(\alpha)}, b_{(\beta)}] c = (a_{(\beta^{-1} \alpha)}
  b)_{(\beta)} c$.

\end{deflist}

Now we observe that one can forget that $\Gamma$ is a
subgroup of $\CC^\times$ and that $R$ is a space of fields.   What
remains is a (not necessarily abelian) group $\Gamma$, its
representation space $R$ and $\CC$-bilinear products
$a_{(\alpha)} b$ on $R$ for each $\alpha \in \Gamma$, such that
the axioms (C0)--(C3) hold.  This is the basic
object of our study, called a \emph{$\Gamma$-conformal algebra}.

Recall that the OPE \req{0.1} of local fields similarly leads to the
notion of a conformal algebra~\cite{K1}.  However, conformal
algebra and $\Gamma$-conformal algebras have dramatically
different properties.

Our first result on $\Gamma$-conformal algebras is Theorem~3.1
which gives their classification in terms of admissible pairs
$(\fg, \varphi)$, where $\fg$ is a Lie algebra and $\varphi$ is a
homomorphism of $\Gamma$ to the group $\Aut \fg$ of automorphisms
of $\fg$ such that for any $a, b \in \fg$ one has
\begin{equation}
  \label{eq:0.6}
  [T_\alpha a, b] = 0 \hbox{~for all but finitely many~}
  \alpha \in \Gamma \, .
\end{equation}

Of course, if $\Gamma$ is a finite subgroup of $\Aut \fg$, then
\req{0.6} holds.
The associated $\Gamma$-conformal algebra is called a twisted
current conformal algebra.  A more interesting example
corresponds to the admissible pair $(g \ell_\infty, \ZZ)$, where
$g \ell_\infty$ is the Lie algebra of all $\ZZ \times \ZZ$-matrices
$(a_{i j})$ over $\CC$ with a finite number of non-zero entries
and the action of $\ZZ$ is by translations along the diagonal:  $n
\cdot (a_{i, j}) = (a_{i+n, j+n})$, $n \in \ZZ$.  This
$\ZZ$-conformal algebra is denoted by $g c_1 (\ZZ)$ for reasons
explained below.

To a $\Gamma$-conformal algebra and a homomorphism $\chi$ of
$\Gamma$ to $\CC^\times$ (or, more generally, to the group of
meromorphic transformations, see \S~\ref{sec:5}) one canonically
associates an infinite-dimensional Lie algebra.  This
construction produces many known (and a lot of new) examples of
infinite-dimensional Lie algebras which appear in deformation of
conformal field theories.  For example, the Lie algebra of
$q$-differential operators on the circle (the sin-algebra) is
obtained from the $\ZZ$-conformal algebra $g c_1 (\ZZ)$ by using
$\chi: \ZZ \to \CC^\times$ defined by $\chi (n) = q^n$.  Taking
the same $\ZZ$-conformal algebra, but the homomorphism $\chi: \ZZ
\to G L_2 (\CC)$ defined by $\chi (n) z = z + n$, $z \in \CC$,
instead, produces the Lie algebra of pseudodifferential operators
on the circle (see \S~\ref{sec:5}).

Representations of $\Gamma$-conformal algebras are classified by
equivariant\break
$(\fg, \Gamma)$-modules satisfying a finiteness
condition (Theorem~4.1).  On the other hand, to any
representation $\pi$ of the group $\Gamma$ we canonically
associate the general $\Gamma$-conformal algebra $g c (\pi,
\Gamma)$, which plays the role of the general linear group in the
sense that representations of a $\Gamma$-conformal algebra $R$
correspond to homomorphisms $R \to g c (\pi, \Gamma)$
(Theorem~4.2).  The most important case is when $\pi$ is a free
$\CC [\Gamma]$-module of rank $N$.  The corresponding general
$\Gamma$-conformal algebra, denoted by $g c_N (\Gamma)$, is a
generalization of $g c_1 (\ZZ)$ mentioned above. It is an analogue
of $g c_N$ in the theory of conformal algebras \cite{K2}.

In the present paper we consider the simplest case of simple
poles and of one indeterminate~$z$.  It is straightforward to
generalize this to the case of multiple poles and the case of
several indeterminates.  This will be discussed elsewhere.

\section{$S$-local formal distributions}
\label{sec:1}

Recall that a formal distribution with coefficients in a vector
space $U$ is an expression of the form
\begin{displaymath}
  a (z) = \sum_{n \in \ZZ} a_n z^{-n - 1} \, ,
  \hbox{~where~}
  a_n \in \ZZ \, .
\end{displaymath}
They form a vector space denoted by $U \left[ [z, z^{-1}]
\right]$.  We shall use the standard notation:
\begin{displaymath}
  \Res_{z = 0} a(z) = a_0 \, .
\end{displaymath}
Similarly, a formal distribution in $z$ and $w$ is an expression
of the form
\begin{displaymath}
  a (z, w) = \sum_{m,n \in \ZZ} a_{m, n}
             z^{-m - 1} w^{-n - 1} \, ,
             \hbox{~where~}
             a_{m, n} \in U \, .
\end{displaymath}
It may be viewed as a formal distribution in $z$ with
coefficients in $U \left[ [w, w^{-1}] \right]$.

Recall that the formal $\delta$-function is the following formal
distribution in $z$ and $w$ with coefficients in $\CC$:
\begin{displaymath}
  \delta (z - w) = \sum_{n \in \ZZ} z^{n - 1} w^{-n} \, .
\end{displaymath}
It has the following two basic properties:
\begin{eqnarray}
  (z - w) \delta (z - w) &=& 0
      \label{eq:1.1} \\
  \Res_{z = 0} a(z) \delta (z - w) &=& a (w)\,
      \hbox{~for~}\, a(z) \in U \left[ [z, z^{-1}] \right] \, .
      \label{eq:1.2}
\end{eqnarray}

We shall use the formal distribution
\begin{displaymath}
  \delta (\alpha z - \beta w)
     = \sum_{n \in \ZZ} (\alpha z)^{n - 1} (\beta w)^{-n} \, ,
     \quad
     \alpha, \beta \in \CC^\times \, ,
\end{displaymath}
and the following obvious property
\begin{equation}
  \label{eq:1.3}
  \delta (\alpha z - \beta w)
     = \alpha^{-1} \delta (z - \alpha^{-1} \beta w)
     = \beta^{-1} \delta (\alpha \beta^{-1} z - w) \, .
\end{equation}
We shall use also the following equality of distribution in
$z_1$, $z_2$ and $z_3$:
\begin{equation}
  \label{eq:1.4}
  \delta (z_1 - z_2) \delta (z_2 - z_3)
     = \delta (z_1 - z_2) \delta (z_1 - z_3)
\end{equation}

Let $S$ be an $N$-element set of (distinct) non-zero complex
numbers.  We shall use the following notation:
\begin{equation}
  \label{eq:1.5}
  (z - w)^N_S = \prod_{\alpha \in S} (z - \alpha w)
\end{equation}

\begin{definition}
  \label{def:1.1}
  Two formal distributions $a (z)$ and $b (z)$ with coefficients
  is Lie algebra $\fg$ are called $S$-\emph{local} if in $\fg
  \left[ [z, z^{-1}, w, w^{-1}] \right]$ one has:
  \begin{equation}
    \label{eq:1.6}
    (z - w)^N_S [a(z), b(w)] = 0 \, .
  \end{equation}
\end{definition}

Note that $S$-local formal distributions are $S'$-local for any
finite subset $S'$ of $\CC^\times$ containing $S$.

\begin{proposition}
  \label{prop:1.1}
  If $a (z)$ and $b(z)$ are $S$-local formal distributions, then
  there exists a unique decomposition of the form
  \begin{equation}
    \label{eq:1.7}
    [a (z), b(w)] = \sum_{\alpha \in S} c_\alpha (w) \delta (z -
    \alpha w) \, .
  \end{equation}
  The formal distributions $c_\alpha (w)$ are given by the
  formula
  \begin{equation}
    \label{eq:1.8}
    c_\alpha (w) = \Res_{z = 0} P_{S, \alpha} \left( \frac{z}{w} \right)
    [a(z), b(w)] \, ,
  \end{equation}
  where
  \begin{equation}
    \label{eq:1.9}
    P_{S, \alpha} (u) = \prod_{\beta \in S \atop \beta \neq \alpha}
                        \frac{u - \beta}{\alpha - \beta}
  \end{equation}
\end{proposition}

Proof of this proposition is immediate by the following lemma.

\begin{lemma}
  \label{lem:1.1}
%%  \item %%%
(a)
    Each formal distribution $a (z, w)$ can be uniquely
    written in the form:
    \begin{eqnarray}
      a (z, w) &=& \sum_{\alpha \in S} c_\alpha (w)
                   \delta (z - \alpha w) + b (z, w)
                   \label{eq:1.10} \\
%%      \noalign{\hbox{\hspace*{25pt} where}}
      \noalign{\hbox{ where}}
      c_\alpha (w) &=& \Res_{z = 0} P_{S, \alpha}
                       \left( \frac{z}{w} \right) a (z, w)
                       \label{eq:1.11} \\
%%      \noalign{\hbox{\hspace*{25pt} and}}
      \noalign{\hbox{ and}}
      \Res_{z = 0} P_{S, \alpha} \left( \frac{z}{w} \right)
                b (z, w) &=& 0
                \quad \hbox{for all~}
                \alpha \in S \, .
                \label{eq:1.12}
    \end{eqnarray}
    %
%%
%%  \item %%%
(b)
    If $b (z, w)$ satisfies \req{1.12} and $(z - w)^N_S b (z, w)
    = 0$\,,\, then $b (z, w) \equiv 0$.
%%  \end{enumerate}
\end{lemma}

\begin{proof}
  Consider the following operator on the space of formal
  distribution in $z$ and $w$:
  \begin{displaymath}
    \pi_S a (z,w) = \sum_{\alpha \in S}
                    c_\alpha (w) \delta (z - \alpha w) \, ,
  \end{displaymath}
  where $c_\alpha (w)$ are given by \req{1.11}.  It is clear
  using \req{1.2} that $\pi^2_S = \pi_S$, which implies (a).
  Furthermore, suppose that $(z - w)^N_S b (z, w) = 0$.  It
  follows from \cite[corollary 2.2]{K1} that for each $\alpha \in S$:
  \begin{displaymath}
    P_{S, \alpha} b(z, w) = d_\alpha (w) \delta (z - \alpha w)
  \end{displaymath}
  for some formal distribution $d_\alpha (w)$.  If also
  \req{1.12} holds, then $\Res_{z=0} d_\alpha (w)$ $\delta (z -
  \alpha w) = 0$, hence $d_\alpha (w) = 0$ by \req{1.2}, and
  \begin{equation}
    \label{eq:1.13}
    P_{S, \alpha} \left( \frac{z}{w} \right) b (z, w) = 0
    \quad \hbox{for all} \quad
    \alpha \in S \, .
  \end{equation}
  This implies that $b (z, w) = 0$ since the polynomials, $P_{S,
    \alpha}$ are relatively prime.  (Indeed $\sum_{\alpha \in S}
  P_{S, \alpha} u_\alpha = 1$ for some polynomials $u_\alpha$.
  Multiplying both sides of \req{1.13} by $u_\alpha \left(
    \frac{z}{w} \right)$ and summing over $\alpha$, we get $b (z,
  w) = 0$.)  This proves~(b).
\end{proof}

\begin{corollary}
  \label{cor:1.1}
  The formal distributions $a (z)$ and $b (z)$ are $S$-local iff
  there exists a decomposition of the form~\req{1.7}.
\end{corollary}

\begin{proof}
  Follows from Proposition~\ref{prop:1.1} and~\req{1.1}.
\end{proof}

\begin{remark}
  \label{rem:1.1}
  The coefficients $c_\alpha (w)$ are independent of the choice
  of $S$ for which \req{1.6} holds.  This follows from the
  uniqueness of the decomposition \req{1.10}.
\end{remark}

\begin{remark}
  \label{rem:1.2}
  Formula \req{1.7}, written out in modes, looks as follows:
  \begin{displaymath}
    [a_m, b_n] = \sum_{\alpha \in S} \alpha^m c_{\alpha, m+n} \, .
  \end{displaymath}
\end{remark}

\section{$\alpha$-products of $\Gamma$-local formal
  distributions}
\label{sec:2}

Let $\Gamma$ be a subgroup of $\CC^\times$.  Two formal distributions
$a (z)$ and $b (z)$ are called $\Gamma$-local if they are
$S$-local for some finite subset $S$ of $\Gamma$.  For $\alpha \in
\Gamma$ define the $\alpha$-product $(a_{(\alpha)} b) (w)$ as the
coefficient $c_\alpha (z)$ in the expansion~\req{1.7}.  In other
words,
\begin{equation}
  \label{eq:2.1}
  \left( a_{(\alpha)} b \right) (w) = \Res_{z=0} P_{S, \alpha}
      \left( \frac{z}{w} \right) [a(z), b(w)] \, .
\end{equation}
Due to Remark~\ref{rem:1.1}, the $\alpha$-product is independent of
the choice of $S$.

For $\beta \in \Gamma$ introduce the following operator:
\begin{displaymath}
  T_\beta (a (z) ) = \beta a (\beta z) \, .
\end{displaymath}
It is clear that $T_\beta$ preserves $\Gamma$-locality.

We collect below the main properties of $\alpha$-products.

\begin{proposition}
  \label{prop:2.1}
  For $\Gamma$-local formal distributions one has the following
  properties $(\alpha, \beta \in \Gamma)$:

    \alphaparenlist
    \begin{enumerate}
    \item %%(a)
      (Translation invariance)
      $(T_\alpha a)_{(\beta)} b = a_{(\alpha \beta)} b$, and
      $a_{(\beta)} T_\alpha b = T_\alpha \left( a_{(\beta
          \alpha^{-1})} b \right)$.

    \item %%(b)
      (Skewsymmetry)
      $a_{(\alpha)} b = -T_\alpha (b_{(\alpha^{-1})} a)$.

    \item %%(c)
      (Jacobi identity)
      $a_{(\alpha)} (b_{(\beta)} c) = a_{(\alpha \beta^{-1})}
      b)_{(\beta)} c + b_{(\beta)} (a_{(\alpha)} c)$ provided
      that all the suitable localities hold.
    \end{enumerate}
\end{proposition}

\begin{proof}
  \alphaparenlist
  \begin{enumerate}
  \item %%(a)
    We have:
    \begin{eqnarray*}
      \left[ T_\alpha a(z), b(w) \right]
         &=& \sum_{\beta \in \Gamma}
             \left( (T_\alpha a)_{(\beta)} b \right)
             (w) \delta (z - \beta w) \, , \\
      \left[ T_\alpha a(z), b(w) \right]
         &=& [ \alpha a (\alpha z), b(w) ]
             = \sum_{\gamma \in \Gamma} \alpha (a_{(\gamma)} b) (w)
             \delta (\alpha z - \gamma w) \\
      &=& \sum_{\gamma \in \Gamma} (a_{(\gamma)} b) (w)
             \delta (z - \alpha^{-1} \gamma w) \, .
    \end{eqnarray*}
    Comparing coefficients, we get the first equation of (a).

  \item %%(b)
    We have using \req{1.7} and \req{1.3}:
    \begin{eqnarray*}
      [ a(z), b(w)] &=& - [b(w), a(z)]
         = - \sum_{\beta \in \Gamma}
           (b_{(\beta)} a) (z) \delta (w - \beta z) \\
       &=& - \sum_{\beta \in \Gamma} \beta (b_{(\beta^{-1})} a)
            (\beta w) \delta (z - \beta w)
    \end{eqnarray*}
    Comparing with \req{1.7}, we obtain (b).  The second equation
    of (a) follows from the first and (b)

  \item %%(c)
    We have:
    \begin{eqnarray*}
      [a(z), [b(w), c(t)] ]
         &=& \sum_{\alpha, \beta \in \Gamma}
             \left( a_{(\alpha)} (b_{(\beta)} c) \right) (t) \delta
             (z - \alpha t) \delta (w - \beta t) \, , \\{}
      [[a(z), b(w)], c(t)]
         &=& \sum_{\alpha, \beta \in \Gamma}
             \left( (a_{(\alpha \beta^{-1})} b)_{(\beta)} c \right) (t)
             \delta (z - \alpha \beta^{-1} w) \delta (w - \beta t) \\{}
         &=& \sum_{\alpha \in \Gamma}
             \left( (a_{(\alpha \beta^{-1})} b)_{(\beta)} c \right) (t)
             \delta (z - \alpha t) \delta (w - \beta t) \, .
    \end{eqnarray*}
    We have used here \req{1.3} and \req{1.4}.  We also have:
    \begin{displaymath}
      [b(w), [a(z), c(t)] ] =
         \sum_{\alpha, \beta \in \Gamma}
         b_{(\beta)} (a_{(\alpha)} c) (t) \delta
         (w - \beta t) \delta (z - \alpha t) \, .
    \end{displaymath}
  \end{enumerate}
  Equating the first equality with the sum of the remaining two
  and using Jacobi identity in $\fg$, we obtain equality (c) due
  to linear independence of $\delta (w - \beta t) \delta (z -
  \alpha t)$ for $\alpha, \beta \in \Gamma$.
\end{proof}

\section{$\Gamma$-conformal algebras and Lie algebras of
  $\Gamma$-local formal distributions}
\label{sec:3}

The above considerations motivate the following definitions.  Let
$\Gamma$ be a group and let $\CC [\Gamma]$ be the group algebra
of $\Gamma$ with basis $T_\alpha$ $(\alpha \in \Gamma)$, so that
$T_\alpha T_\beta = T_{\alpha \beta}$, $T_1 = 1$.  Since $\Gamma$
is not necessarily abelian we must distinguish the case of left
and right $\CC [\Gamma]$-modules.

\begin{definition}
  \label{def:3.1}
  A left $\CC [\Gamma]$-module $R$ is called a left
  \emph{$\Gamma$-conformal algebra} if it is equipped with a
  $\CC$-bilinear product $a_{(\alpha)} b$ for each $\alpha \in
  \Gamma$ such that the following axioms hold $(a, b, c \in R \,
  , \; \alpha, \beta \in \Gamma)$:

\begin{deflist}{C}

\item %% C0
  $a_{(\alpha)} b = 0$ for all but finitely many $a \in \Gamma$,

\item %% C1
  $(T_\alpha a)_{(\beta)} b = a_{(\beta \alpha)} b$,

\item %% C2
  $a_{(\alpha)} b = - T_\alpha (b_{(\alpha^{-1})} a)$,

\item %% C3
  $a_{(\alpha)} (b_{(\beta)} c) = (a_{(\beta^{-1} \alpha)} b)_{(\beta)} c +
  b_{(\beta)} (a_{(\alpha)} c)$.
\end{deflist}

\end{definition}

A right $\CC [\Gamma]$-module $R$ is called a \emph{right
  $\Gamma$-conformal algebra} if axioms (C1), (C2) and (C3) are
replaced by:

\begin{optenum}{(C2)$_{\mathrm R}$}

\item [(C1)$_{\mathrm R}$]
  $(a T_\alpha)_{(\beta)} b = a_{(\alpha \beta)} b$,

\item [(C2)$_{\mathrm R}$]
  $a_{(\alpha)} b = - (b_{(\alpha^{-1})} a) T_\alpha$,

\item [(C3)$_{\mathrm R}$]
  $a_{(\alpha)} (b_{(\beta)} c) = (a_{(\alpha \beta^{-1})}
  b)_{(\beta)} c + b_{(\beta)} (a_{(\alpha)} c)$.

\end{optenum}

Unless otherwise specified, we will consider left $\Gamma$-conformal
algebras and will drop the adjective left.

\begin{remark}
  \label{rem:3.1}
    \alphaparenlist
    \begin{enumerate}
    \item %%(a)
      Properties (C1) and (C2) (resp. (C1)$_{\mathrm R}$ and
      (C2)$_{\mathrm R}$) imply

      \begin{optenum}{(C1)$_{\mathrm R}$}

      \item [(C1$'$)] \hspace{-1ex}
        $a_{(\beta)} T_\alpha b = T_\alpha (a_{(\alpha^{-1} \beta)} b)$
        (resp. (C1$'$)$_{\mathrm R}$
        $a_{(\beta)} (b T_\alpha) = (a_{(\beta \alpha^{-1})} b) T_\alpha)$

     \end{optenum}

   \item %%(b)
     (C1) and (C1$'$) (resp. (C1)$_{\mathrm R}$ and
     (C1$'$)$_{\mathrm R}$) imply
     \begin{displaymath}
       (T_\alpha a)_{(\beta)} (T_\alpha b)
          = T_\alpha (a_{(\alpha^{-1} \beta \alpha)} b)
            \quad
            (\resp (a T_\alpha)_{(\beta)} (b T_\alpha) =
            (a_{(\alpha \beta \alpha^{-1})} b) T_\alpha) \, .
     \end{displaymath}
     In particular, each $T_\alpha$ is an automorphism of the
     product $a_{(1)} b$ for each $\alpha \in \Gamma$.

   \item %%(c)
     With respect to the product $a_{(1)} b$, $R$ is a Lie
     algebra over $\CC$.  Due to (C1) all other products are
     expressed via the $1$-product:
     \begin{displaymath}
       a_{(\alpha)} b = (T_\alpha a)_{(1)} b
     \end{displaymath}
    \end{enumerate}
\end{remark}

Due to Remark~\ref{rem:3.1}, we can associate to a
$\Gamma$-conformal algebra $R$ a pair $(\fg, \varphi)$, where
$\fg$ is a Lie algebra over $\CC$ whose underlying space is $R$
and bracket is $[a, b] = a_{(1)} b$, and $\varphi$ is a
homomorphism of $\Gamma$ to the group $\Aut (\fg)$ such that
\begin{equation}
  \label{eq:3.1}
  [T_\alpha a, b] = 0
  \quad
  \hbox{for all but finitely many $\alpha \in \Gamma$.}
\end{equation}
We call such $(\fg, \varphi)$ an \emph{admissible pair}.
Conversely, given an admissible pair $(\fg, \varphi)$ we denote
by $R (\fg, \varphi)$ the $\CC [\Gamma]$-module $\fg$ with
$\CC$-bilinear products
\begin{equation}
  \label{eq:3.2}
  a_{(\alpha)} b = [T_\alpha a, b] \
\end{equation}
It is straightforward to check that $R (\fg, \varphi)$ is a
$\Gamma$-conformal algebra.

We have proved the following result:

\begin{theorem}
  \label{theor:3.1}
  $\Gamma$-conformal algebras are classified by admissible pairs
  $(\fg, \varphi)$.
\end{theorem}

Fix now a homomorphism $\chi: \Gamma \to \CC^\times$.  Let $\fs$
be a Lie algebra and let $R$ be a family of pairwise $\chi
(\Gamma)$-local formal distributions with coefficients in $\fs$,
which span $\fs$ over $\CC$.  Assume that $R$ is stable under all
$T_\alpha$ and all $\alpha$-products for $a \in \chi (\Gamma)$.
Then $\fs$ is called a \emph{Lie algebra} of \emph{$\chi
  (\Gamma)$-local formal distributions}.

It follows from Proposition~\ref{prop:2.1} that $R$ is a $\chi
(\Gamma)$-conformal algebra.  Conversely, given a
$\Gamma$-conformal algebra $R$, we associate to it a Lie algebra
$\fs (R)$ of $\chi (\Gamma)$-local formal distributions as
follows.  Consider a vector space over $\CC$ with the basis
$a_n$, where $a \in R$ and $n \in \ZZ$, and denote by $\fs (R,
\chi)$ the quotient of this space by the $\CC$-span of elements
of the form $(n \in \ZZ)$:
\begin{mathletters}
  \label{eq:3.3}
  \begin{equation}
    (\lambda a + \mu b)_n - \lambda a_n - \mu b_n \, ,
       \quad
       \hbox{where }
       \lambda, \mu \in \CC \, ,
       \; a, b \in R \, ,
       \label{eq:3.3a} %%\\
  \end{equation}
  \begin{equation}
    (T_\alpha a)_n - \chi (\alpha)^{-n} a_n \, ,
       \quad
       \hbox{where }
       \alpha \in \Gamma \, ,
       \; a \in R \, .
    \label{eq:3.3b}
  \end{equation}
\end{mathletters}
Then it is straightforward to check that the formula (cf.\
Remark~\ref{rem:1.2})
\begin{equation}
  \label{eq:3.4}
  [a_m, b_n] = \sum_{\gamma \in \Gamma} \chi (\gamma)^m
  (a_{(\gamma)} b)_{m+n}
\end{equation}
gives a well-defined structure of a Lie algebra on $\fs (R,
\chi)$ of $\chi (\Gamma)$-local formal distributions $a (z) = \sum_{n
  \in \ZZ} a_n z^{-n-1}$, $a \in R$.
The associated $\chi (\Gamma)$-conformal algebra is denoted by $R
(\fs)$.

The relations between the constructed objects can be summarized
by the following diagram:
  \begin{center}
      \setlength{\unitlength}{.5ex}
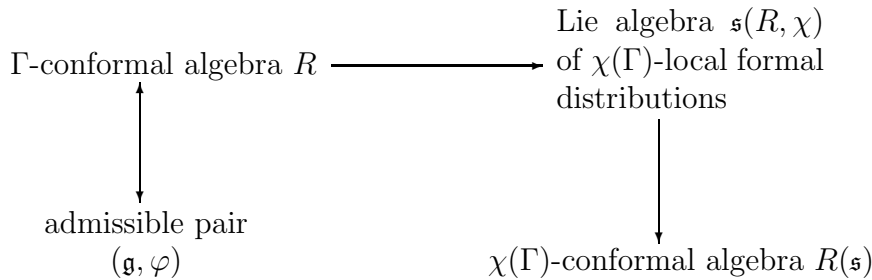
      \begin{picture}(100,45)(0,5)
        \put(-10,40){\makebox(0,0)[l]{$\Gamma$-conformal
            algebra $R$}}
        \put(70,40){\parbox[l]{1.4in}{Lie algebra $\fs (R, \chi)$
            of $\chi (\Gamma)$-local formal distributions}}
        \put(37,40){\vector(1,0){31}}
        \put(85,32){\vector(0,-1){18}}
        \put(-5,10){\shortstack[c]{admissible pair \\
            $(\fg, \varphi)$ }}
        \put(60,10){\parbox[tr]{30ex}{$\chi (\Gamma)$-conformal
            algebra $R (\fs)$}}
        \put(9,20.5){\vector(0,1){17}}
        \put(9,37){\vector(0,-1){17}}
      \end{picture}
\end{center}

Now we turn to the most important examples of $\Gamma$-conformal
algebras and related constructions.

\begin{example}
\rm
  \label{ex:3.1}
  Let $(\fg, \ZZ / N \ZZ)$ be a finite-dimensional (simple) Lie
  algebra $\fg$ with the action of the group $\Gamma = \ZZ / N
  \ZZ$ by automorphisms of $\fg$.  Due to this action $T_\alpha$,
  $\alpha \in \Gamma$, we have the eigenspace decomposition
  ($\hat{\Gamma}$ is the group of characters of $\Gamma$):
  \begin{displaymath}
    \fg = \subsym{\oplus}{j \in \hat{\Gamma}} \fg_j \, .
  \end{displaymath}
The $\Gamma$-conformal algebra $R (\fg, \varphi)$ is defined as
$\CC [\Gamma]$-module with underlying space $\fg$ and
$\CC$-bilinear products $a_{(\alpha)} b = [T_\alpha a, b]$, $a,
b \in \fg$, $\alpha \in \ZZ / N \ZZ$.  Fix a homomorphism $\chi :
\Gamma \to \CC^\times$ such that $\chi (n) = \epsilon^n$, where
$\epsilon \in \CC$ is an $N$-th root of $1$.  The corresponding
Lie algebra of $\chi (\Gamma)$-local formal distributions $\fs
(R, \chi)$ is nothing but
the twisted affine algebra $[ K ]$ with the commutation relations:
\begin{displaymath}
  [a_m, b_n] = \sum_{k \in \ZZ / N \ZZ} \epsilon^{k m} [T_k a,
  b]_{m + n}
\end{displaymath}
Equivalently:
\begin{displaymath}
[a(z), b(w)] = \sum_{\alpha \in \ZZ / N \ZZ} [T_\alpha a, b] (w)
\delta (z - \alpha w) \, .
\end{displaymath}
Suppose that $a \in \fg_j$, then $T_k a = \epsilon^{k j} a$ and
\begin{displaymath}
  [a_m, b_n] = \left\{
      \begin{array}{c@{\quad}l}
        0 & \hbox{if } m \neq -j \bmod N \\
        N [a, b]_{m+n} & \hbox{if } m = -j \bmod N
      \end{array}
    \right.
\end{displaymath}
Thus we put $a_m = 0$ for $m \neq -j \bmod N$, which is the same
as to take a quotient by relations \req{3.3b}.
\end{example}

\begin{example}
\rm
  \label{ex:3.2}
  Consider $\Gamma = \ZZ$ and the admissible pair $(\fg,
  \varphi)$, where $\fg = g \ell_\infty$ is the Lie algebra of all $
 \ZZ \times \ZZ$ matrices over $\CC$ with finitely many non-zero
 entries.  Let $E_{i j}$ $(i,j \in \ZZ)$  be its standard basis,
 with the usual relations:
 \begin{displaymath}
   [E_{i j}, E_{k \ell}] =
       \delta_{j, k} E_{i \ell}- \delta_{\ell, i} E_{k j} \, .
 \end{displaymath}
 Let $T$ be the image of $1 \in \ZZ$ under $\varphi$.  Consider
 two different actions of $T$ on $g \ell_\infty$.

 \alphalist
 \begin{enumerate}
 \item %%a)
   $T : E_{i j} \to E_{i + 1, j + 1}$. \\
   In this case $g \ell_\infty$ is a free $\CC [T, T^{-1}]$-module
   with the generators $A^m = E_{0, m}$, $m \in \ZZ$.  As follows
   from Definition~\ref{def:3.1} of a $\Gamma$-conformal algebra,
   it is sufficient to define the $\Gamma$-products only on
   generators of the $\CC [\Gamma]$-module.  So, for the case
   under consideration we have:
   \begin{eqnarray}
     \label{eq:3.7}
     A^m_{(r)} A^n &=& [T^r A^m, A^n]
          = [E_{r, m + r}, E_{o, n}] \nonumber \\
        &=& \delta_{r, -m} T^{-m} A^{m + n}
           - \delta_{r, n} A^{m + n} \, ,
           \quad
           r \in \ZZ
   \end{eqnarray}
   Take a homomorphism $\chi = \chi_q : \ZZ \to \CC^\times$
   defined by $\chi_q (1) = q \in \CC^\times$.  The corresponding
   Lie algebra $\fs (R, \chi)$ is the Lie algebra with the basis
   $A^m_k$, $m, k \in \ZZ$, and the commutation relations:
   \begin{equation}
     \label{eq:3.8}
     [A^m_k, A^n_\ell] =
        \sum_{r \in \ZZ} q^{r k} (A^m_{(r)} A^n)_{k + \ell}
        = (q^{m \ell} - q^{k n}) A^{m + n}_{k + \ell} \, .
   \end{equation}
   This is the Lie algebra of $q$-pseudodifferential operators on
   the circle (sin-algebra).

   \item %% b)
     Now define the action of $\varphi (1) = \tilde{T}$ by
     $\tilde{T} (E_{i, j}) = - E_{j + 1, i + 1}$.  We have
     $\tilde{T} = \epsilon T$, where $\epsilon$ is an order $2$
     automorphism of $g \ell_\infty$ defined by $\epsilon (E_{i
       j}) = - E_{j i}$.  In this case $g \ell_\infty$ is a
     free $\CC [\tilde{T}, \tilde{T}^{-1}]$-module with
     generators $A^m = E_{0, m}$ and $B^m = \epsilon
     (A^m) = - E_{m, 0}$, $m \in \ZZ_+$.  Let us compute the
     $r$-product, $r \in \ZZ$, for the elements of the basis of $
     \CC [\Gamma]-module$.  Suppose $n - m \ge 0$ (for $n - m < 0$
     we can use axiom (C2)):
     \begin{eqnarray*}
       A^m_{(r)} A^n &=& [\tilde{T}^r A^m, A^n]
             \nonumber \\
        &=&
          \left\{
            \begin{array}{l@{\quad}l}
              \delta_{r, -m} T^{-m} A^{m + n} - \delta_{r, m}
                 A^{m + n}\, ,
                 & m, n \hbox{ even;} \\
              \delta_{r, -m} T^{-m} A^{m + n} + \delta_{r, n - m}
                 A^{n - m}\, ,
                 & m \hbox{ even, } n \hbox{ odd;} \\
              -\delta_{r, n} A^{m + n} + \delta_{r, n - m}
                 A^{n - m}\, ,
                 & m \hbox{ odd, } n \hbox{ even} \\
              0 &
                m, n \hbox{ odd}. \\
            \end{array}
          \right. \nonumber \\
        B^m_{(r)} B^n &=& \epsilon [\tilde{T}^r A^m, A^n]
          \label{eq:3.7b} \\
        B^m_{(r)} A^n &=&
          \left\{
            \begin{array}{l@{\quad}l}
              -\delta_{r, 0} \tilde{T}^{m} A^{n -m} + \delta_{r, n-m}
                 A^{n-m}\, ,
                 & m, n \hbox{ even;} \\
              -\delta_{r, 0} \tilde{T}^{m} A^{n -m} - \delta_{r, n}
                 A^{m+n}\, ,
                 & m \hbox{ even, } n \hbox{ odd;} \\
              -\delta_{r, 0} \tilde{T}^{m} B^{n-m} + \delta_{r, -m}
                 \tilde{T}^{-m} B^{n+m}\, ,
                 & m \hbox{ odd, } n \hbox{ even;} \\
              -\delta_{r, 0} \tilde{T}^{n} B^{n-m} + \delta_{r, n-m}
                 A^{n+m} \\
                 \quad {} + \delta_{r, -m}
                 \tilde{T}^{-m} B^{n+m} - \delta_{r, n}
                 A^{m + n} \, ,
                 & m,n \hbox{ odd}.
            \end{array}
          \right.
     \end{eqnarray*}
     The corresponding Lie algebra $\fs (R, \chi_q)$ is a Lie
     algebra with basis $\{ A^m_k \}$ and $\{ B^n_\ell \}$, $m,n \in
     \ZZ_+$, $k, \ell \in \ZZ$.  From \req{3.4} and \req{3.7} we
     have for $n - m \ge 0$:
     \begin{eqnarray*}
       [A^m_k, A^n_\ell] &=&
          \left\{
            \begin{array}{l@{\quad}l}
              (q^{m \ell} - q^{n k}) A^{m+n}_{k+\ell}\, ,
                 & m,n \hbox{ even}, \\
              q^{m \ell} A^{m+n}_{k+\ell} + q^{k (n-m)}
                 A^{n-m}_{k+\ell}\, ,
                 & m \hbox{ even, } n \hbox{ odd;} \\
              -q^{n k} A^{m+n}_{k+\ell} + q^{k (n-m)}
                 A^{n-m}_{k+\ell}\, ,
                 & m \hbox{ odd, } n \hbox{ even;} \\
              0  & m,n \hbox{ odd,}
            \end{array}
          \right. \\{}
       [B^m_k, B^n_\ell] &=&
          \hbox{the same as for }
          [A^m_n, A^n_\ell]
          \hbox{ with change }
          A \to B \, . \\{}
       [B^m_k, A^n_\ell] &=&
          \left\{
            \begin{array}{l@{\quad}l}
              q^{-m k} (q^{k n} - q^{-m \ell}) A^{n-m}_{k+\ell}\, ,
                 & m,n \hbox{ even}, \\
              -q^{-m (k+\ell)} A^{n-m}_{k+\ell} -q^{k n} A^{m+n}_{k+\ell}\, ,
                 & m \hbox{ even, } n \hbox{ odd;} \\
              -q^{-m (k+\ell)} B^{n-m}_{k+\ell} +q^{m \ell}
                   B^{n+m}_{k+\ell} \, ,
                 & m \hbox{ odd, } n \hbox{ even;} \\
              -q^{-m (k+\ell)} B^{n-m}_{k+\ell} +q^{k (n-m)}
                   A^{n-m}_{k+\ell} \\
                \quad {} +q^{m \ell}
                   B^{n+m}_{k+\ell} -q^{n k}
                   A^{n+m}_{k+\ell}\, ,
                 & m,n \hbox{ odd}.
            \end{array}
          \right.
     \end{eqnarray*}
     This Lie algebra is apparently new.  It contains as a
     subalgebra the sin-algebra $\{ \tilde{A}^k_n \}$, where
     $\tilde{A}^n_k = A^n_k$, $n > 0$, $n$ even; $A^n_k = q^{n
       k} B^{-n}_k$, $n < 0$, $n$ even; and the subalgebra
     isomorphic to $B$-series of Lie algebras of quantum torus
     (see \cite{G-K-L}) for $n$ odd.

\end{enumerate}
\end{example}

\begin{example}             
\rm
  \label{ex:3.3}
  Consider the group $\Gamma = \ZZ_2 \times \ZZ$ and admissible
  pair $(g \ell_\infty, \varphi)$, where homomorphism $\varphi:
  \Gamma \to \Aut g \ell_\infty$ is determined by two elements
  $\epsilon$ and $T$ as in Example~\ref{ex:3.2}~(b).  With
  respect to the $\Gamma$-action $g \ell_\infty$ is a free
  module with the generators $B^m = E_{0, m}$, $m \ge 0$.  The
  corresponding conformal algebra $R (g \ell_\infty, \varphi)$
  is defined by the $\Gamma$-products:
  \begin{eqnarray}
    B^m_{(r)} B^n = [T^r B^m, B^n] &=&
       \delta_{r, -m} T^{-m} B^{n+m} - \delta_{n,r} B^{m+r} \, ,
       \nonumber \\[-1ex]
    \label{eq:3.8b} \\[-0.25ex]
    B^m_{(\epsilon r)} B^n = [\epsilon T^r B^m, B^n] &=&
       \left\{
         \begin{array}{l@{\quad}l}
           -\delta_{r,0} T^m B^{n-m} + \delta_{r, n-m} B^{n-m}\, ,
              & n - m \ge 0 \\
           \delta_{r,0} \epsilon T^n B^{n-m} - \delta_{r, n-m} B^{m-n}\, ,
              & n - m < 0 \, .
         \end{array}
       \right.
     \nonumber
  \end{eqnarray}
  Fix the homomorphism $\chi : \Gamma \to \CC^\times$ such that $
  \chi (\epsilon) = - 1$, $\chi (1) = q \in \CC^\times$.  The Lie
  algebra $\fs (R (g \ell_\infty, \varphi), \chi)$ is a Lie
  algebra with the basis $B^m_k$, $m \in \ZZ_+$, $k \in \ZZ$, and
  commutation relations:
  \begin{eqnarray*}
    [B^m_k, B^n_\ell] &=& (q^{m \ell} - q^{n k}) B^{m + n}_{k + \ell}
       \\
       && \quad {} - \left\{
         \begin{array}{l@{\quad}l}
           (-1)^k (q^{k (n-m)} - q^{-m (k+\ell)}) B^{n-m}_{k+\ell} \, ,
              & n-m > 0 \, , \\
           (-1)^k ((-q)^{-n (k+\ell)} \\
             \quad {} - (-1)^{-n (k+\ell)}
              (-q)^{m (k+\ell)}) B^{m-n}_{k+\ell}\, ,
              & n-m < 0
         \end{array}
       \right.
  \end{eqnarray*}
  After the renormalization of the basis $\tilde{B}^m_k =
  q^{\frac{m k}{2}} B^M_k$ we get exactly the commutation
  relations for $B$-series of sin Lie algebras, introduced in
  \cite{G-K-L}:
  \begin{displaymath}
    [\tilde{B}^m_k, \tilde{B}^n_\ell]
       = \left (q^{\frac{m \ell - n k}{2}}
           - q^{\frac{- (m \ell - n k)}{2}}
         \right)
         \tilde{B}^{n+m}_{k+\ell}
         - (-1)^k
         \left( q^{\frac{k n + m \ell}{2}}
           - q^{-\frac{m \ell + k n}{2}}\right) \tilde{B}^{n-m}_{k+\ell}
\end{displaymath}
\end{example}
%%%%%%%%%%%%%%%%%%%%%%%%%%%%%%%%%%%%%%%%%%%%%%%%%%%%%%%%%%%%%%%%%%
%%%%%%%%%%%%%%%%%%%%%%%%%%%%    Jan's part     %%%%%%%%%%%%%%%%%%%%%%%%%%%%
%%%%%%%%%%%%%%%%%%%%%%%%%%%%%%%%%%%%%%%%%%%%%%%%%%%%%%%%%%%%%%%%%%%%%%%%%%%

\goodbreak
\section  {Representations of $\Gamma$-conformal algebras and the
general $\Gamma$-conformal algebra $gc (\pi, \Gamma)$.}
\label{sec:4}

\begin{definition}
  \label{def:4.1}
  A (left) \emph{module} over a $\Gamma$-conformal algebra $R$ is
  a $\CC [\Gamma]$-module $M$ with a $\CC$-linear map $a \mapsto
  a^M_{\alpha}$ of $R$ to $End_{\CC} M$ for each $\alpha \in
  \Gamma$ such that the following properties hold (here $a, b \in
  R$, $\alpha, \beta \in \Gamma$ and the action of $\Gamma$ on
  $M$ is denoted by $\alpha \mapsto T^M_{\alpha}$):

  \begin{optenum}{(M0)}
    \item[(M0)]
      $a^M_{(\alpha)} v=0 \hbox{ for } v \in M \hbox{ and all
      but finitely many } \alpha$.
    \item[(M1)]
    $(T_{\alpha} a)^M_{(\beta)} = a^M_{\beta \alpha}  \, , \,
    a^M_{(\beta)} T^M_{\alpha} = T^M_{\alpha}
    a^M_{(\alpha^{-1}\beta)}$,
    \item[(M2)]
    $\left[ a^M_{(\alpha)}, b^M_{(\beta)} \right]
     = \left(a_{(\beta^{-1} \alpha)} b \right)^M_{(\beta)}$.
  \end{optenum}
\end{definition}

For example, $a \mapsto a^R_{(\alpha)}$ is the \emph{adjoint
  representation} of $R$ on itself.

It follows from (M1) that
\begin{equation}
  \label{eq:4.1}
  (T_{\alpha} a)^M_{(\beta)}(T^M_{\alpha} v)
  = T^M_{\alpha}(a^M_{(\alpha^{-1} \beta \alpha)} v) \, ,
  \quad
  v \in M \, .
\end{equation}
This leads us to

\begin{theorem}
  \label{theorem:4.1}
  Let $R$ be a $\Gamma$-conformal algebra and let $(\fg, \varphi)$
  be the associated admissible pair.  Then $R$-modules $M$ are
  classified by equivariant $(\fg, \Gamma)$-modules such that for
  any $a \in \fg$ and $v \in M$ one has:
  \begin{equation}
    \label{eq:4.2}
    (T^{\fg}_{\alpha}a) v = 0\, \hbox{ for all but finitely many }
    \alpha \in  \Gamma \, .
  \end{equation}

\end{theorem}

\begin{proof}
  It follows from (M2) that $a \mapsto
  a_{(1)}$ is a representation of $\fg$.  It satisfies
  (\ref{eq:4.1}) due to (M2), and it is equivariant,
  i.e. $(T^{\fg}_{\alpha} a)(T^M_{\alpha} v)=T^M_{\alpha}(av)$ due
  to (\ref{eq:4.1}).  Conversely, given an equivariant $(\fg,
  \Gamma)$-module $M$, it becomes an $R$-module by letting
  \begin{displaymath}
    a^M_{\alpha} =(T^{\fg}_{\alpha} a)^M_{(1)} \hbox{ on } M \, .
  \end{displaymath}
\end{proof}

Consider a representation $\alpha \mapsto T_{\alpha}$ of the
group $\Gamma$ in a vector space $V$ over $\CC$.  We define a
$\Gamma$-\emph{conformal endomorphism} of $V$ as a collection $a=
\left\{ a_{(\beta)} \right\}_{\beta \in \Gamma}$ of
$\CC$-endomorphisms of $V$ such that
\begin{equation}
  \label{eq:4.3}
  a_{(\beta)} T_{\alpha} = T_{\alpha} a_{(\alpha^{-1} \beta)} \, ,
  \quad
  \alpha, \beta \in \Gamma \, ,
\end{equation}
and for each $v \in V$
\begin{equation}
  \label{eq:4.4}
  a_{(\alpha)}v=0 \hbox{ for all but finitely many }\alpha \in
    \Gamma \, .
\end{equation}
We denote the set of all $\Gamma$-conformal endomorphisms of $V$
by $gc (V, \Gamma)$.

Define a $\CC [\Gamma]$-module structure on the space $gc (V,
\Gamma)$ by
\begin{eqnarray}
  \label{eq:4.5}
  (T_{\alpha}a)_{(\beta)} = a_{(\beta \alpha)} \, , \alpha \, ,
  \beta \in \Gamma \, ,
\end{eqnarray}
and $\alpha$-product for each $\alpha \in \Gamma$ by
\begin{equation}
  \label{eq:4.6}
  (a_{(\alpha)} b)_{(\beta)} = [a_{(\beta \alpha)} \, ,
  b_{(\beta)}] \, .
\end{equation}
It is immediate to check that axioms (C1), (C2) and (C3)
of a $\Gamma$-conformal algebra hold (though axiom (C0)
probably doesn't hold in general).

It is clear that we have by definition
\begin{proposition}
\label{proposition:4.1}
To give a $\Gamma$-module $V$ a structure of a module over a
$\Gamma$-conformal algebra $R$ is the same as to give a
homomorphism $R \to gc (V, \Gamma)$ (i.e. a $\CC[\Gamma ]$-module
homomorphism preserving all $\alpha$-products).
 \end{proposition}
We will show now that axiom (C0) holds, hence $gc(V, \Gamma)$ is a
$\Gamma$-conformal algebra, provided that the $\CC[ \Gamma
]$-module $V$ is finitely generated.  Towards this end we shall
give a different construction of $gc(V, \Gamma)$ in the case when
$V$ is a free $\CC [\Gamma ]$-module of rank $1$.

Let $gc_1 (\Gamma)= \bigoplus_{r \in \Gamma} \CC [\Gamma]a^r$ be
a free $\CC[\Gamma]$-module with free generators $a^r$ labeled by
elements of $\Gamma$.  For each $\alpha \in \Gamma$ define the
$\alpha$-product by
\begin{equation}
  \label{eq:4.7}
  a^r {}_{(\alpha)} a^s = \delta_{\alpha,r^{-1}} T_{r^{-1}}
  a^{rs}- \delta_{\alpha, s} a^{sr}\, ,
\end{equation}
extending to the whole $R$ by axioms (C1) and (C1').

\begin{theorem}
  \label{theorem:4.2}
  \alphaparenlist
  \begin{enumerate}
  \item %%(a)
    The $\CC(\Gamma]$-module $gc_1(\Gamma)$ with products
    (\ref{eq:4.7}) is a $\Gamma$-conformal algebra.

  \item  %%(b)
    Let $V= \CC[ \Gamma]\, v$ be a free $\CC[ \Gamma]$-module
    of rank $1$.  Define the action of $gc_1( \Gamma)$ on $V$ by
    letting
    \begin{equation}
      \label{eq:4.8}
      (a^s)^{V}_{(\alpha)} v = \delta_{\alpha,s^{-1}} T_{\alpha}v
    \end{equation}
    and extending to $V$ by (M1).  This gives $V$ a structure
    of a $gc_1( \Gamma)$-module.

  \item  %%(c)
    The $gc_1(\Gamma)$-module $V = \CC [\Gamma] v$ contains all
    $\Gamma$-conformal endomorphisms of $V$.  In particular, the
    $\Gamma$-conformal algebras $gc(V, \Gamma)$ and $gc_1(\Gamma)$
    are isomorphic.

  \end{enumerate}
\end{theorem}

\begin{proof}
  To prove (a) we need to check that the axioms (C0)--(C3) are
  satisfied.  (C0) is obvious.  (C1) and (C2) come from the
  definition (\ref{eq:4.7}) of $gc_1 (\Gamma)$.  Check the skew
  symmetry (C2):
  \begin{eqnarray*}
    ( a^r_{(\alpha)}a^s )_{(\beta)}
        &=& \delta_{\alpha, r^{-1}} (T_{r^{-1}} a^{rs})_{(\beta)}
            - \delta_{\alpha , s} a^{sr}_{(\beta)} \, , \\
    T_{\alpha}(a^s_{(\alpha^{-1})}a^r)_{\beta}
        &=& \delta_{\alpha , s} T_{\alpha s^{-1}} a^{sr}
               - \delta_{\alpha , r^{-1}} T_{r^{-1}} a^{rs}
            = -( \delta_{\alpha , r^{-1}} T_{r^{-1}} a^{rs}
               - \delta_{\alpha , s} a^{sr}) \, .
  \end{eqnarray*}

  For (C3) we need to check that both sides of
  \begin{displaymath}
    a^m_{(\beta \alpha)}(a^r_{(\beta)} a^s)
    = (a^m_{(\alpha)} a^r)_{(\beta)} a^s + a^r_{(\beta)}
    (a^m_{(\beta \alpha)} a^s)
  \end{displaymath}
  are equal.
  \begin{eqnarray*}
    \hbox{LHS:\quad}
    a^m_{(\beta \alpha)}(a^r_{(\beta)}a^s)
       &=&  a^m_{(\beta \alpha)}
           (\delta_{\beta , r^{-1}} T_{r^{-1}} T_{r^{-1}}
           - \delta_{\beta , s} a^{sr}) \\
       &=& \delta_{\beta , r^{-1}}T_{r^{-1}} a^m_{(\alpha)} a^{rs}
           - \delta_{\beta , s} a^m_{(\beta \alpha)} a^{sr} \\
       &=& \delta_{\beta , r^{-1}} \delta_{\alpha , m^{-1}}
           T_{(mr)^{-1}} a^{mrs}
           - \delta_{\beta , r^{-1}} \delta_{\alpha , rs}
           T_{r^{-1}} a^{rsm} \\
       && \quad {} -
          \delta_{\beta , s}
          \delta_{\beta \alpha , m^{-1}} T_{m^{-1}} a^{msr}
          + \delta_{\beta , s} \delta_{\beta \alpha , sr} a^{srm}
  \end{eqnarray*}
  \begin{eqnarray*}
    \lefteqn{\hspace{-4.3em}
      \hbox{RHS:\quad}
      (a^m_{(\alpha)} a^r)_{(\beta)} a^s+a^r_{(\beta)}
      (a^m_{(\beta \alpha)}a^s)
      } \\
    \qquad\qquad&=& \delta_{\alpha , m^{-1}}
        \delta_{\beta \alpha , (mr)^{-1}} T_{(mr)^{-1}} a^{mrs}
        - \delta_{\alpha_1 m^{-1}} \delta_{\beta \alpha , s}
        a^{smr} \\
    && \quad {} -
        \delta_{\alpha , r} \delta_{\beta , (rm)^{-1}} a^{rms}
        + \delta_{\alpha , r} \delta_{\beta , s} a^{srm} \\
    && \quad {}
        + \delta_{\beta \alpha , m^{-1}}
        \delta_{\beta ,r} T_{(rm)^{-1}} a^{rms}
        - \delta_{\beta \alpha ,m^{-1}} \delta_{\alpha^{-1} , ms}
        T_{m^{-1}} a^{msr} \\
    && \quad {}
        - \delta_{\beta \alpha ,s}
        \delta_{\beta ,r^{-1}} a^{rsm} + \delta_{\beta \alpha ,s}
        \delta_{\beta ,sm} a^{smr} \, ,
\end{eqnarray*}
and we have LHS $=$ RHS.

To prove (b) we need to check only axiom (M2) on $V$:
\begin{displaymath}
  (a^r_{(\alpha)} a^s)_{(\beta)} v = [a^r_{(\beta \alpha)},
     a^s_{(\beta)}] v \, .
\end{displaymath}
\begin{eqnarray*}
    \hspace{-4em}
      \hbox{LHS:\quad}
      (a^r_{(\alpha)}a^s)_{(\beta)}v
        &=& ( \delta_{\alpha ,r^{-1}}
            a^{rs}_{\beta \alpha}
            - \delta_{\alpha ,s} a^{sr}_{(\beta)}) v
             \\
        &=& (\delta_{\alpha ,r^{-1}} \delta_{\beta ,s^{-1}}
            T_{(rs)^{-1}}- \delta_{\alpha ,s} \delta_{\beta(s2)^{-1}}
            T_{\beta}) v
\end{eqnarray*}
\begin{eqnarray*}
    \hspace{-1.5em}
      \hbox{RHS:\quad}
      [a^r_{(\beta \alpha)} , a^s_{(\beta)}] v
      &=& ( \delta_{\beta ,s^{-1}} T_{s^{-1}} a^r_{(\alpha)}
          - \delta_{\beta \alpha ,r^{-1}} T_{r^{-1}}
          a^s_{(\alpha^{-1})}) v \\
      &=& \left( \delta_{\beta, s^{-1}} \delta_{\alpha, r^{-1}}
          T_{(r s)^{-1}}
          - \delta_{\beta \alpha, r^{-1}} \delta_{\alpha, s}
          T_{(s r)^{-1}} \right) v
\end{eqnarray*}
so LHS $=$ RHS.  Finally, (c) is clear.
\end{proof}

\begin{corollary}
  \label{cor:4.1}  If $V$ is a finitely generated $\CC[
  \Gamma]$-module, then $gc(V, \Gamma)$ is a conformal algebra.
\end{corollary}

\begin{proof}
  We may assume that $V=\CC [\Gamma ]^{\oplus N}$ is a free
  $\CC[\Gamma]$-module of rank $N$.  Then $gc(V, \Gamma)$ may be
  viewed as a $\CC[\Gamma]$-submodule of $gc_1 (\Gamma^N)$ using
  the diagonal homomorphism $\CC [\Gamma ] \to \CC [\Gamma
  ]^{\oplus N}$.  Corollary now follows from Theorem \ref{theorem:4.2}a.
\end{proof}

If $V$ is a free $\CC [\Gamma ]$-module of rank $N$, we use
notation $gc_N (\Gamma) = gc(V, \Gamma)$.

\begin{example}
\rm
  \label{ex:4.1}
  Consider $\Gamma = \ZZ^N$.  For the $gc_1(\ZZ^N)$ defined by
  (\ref{eq:4.7}) the corresponding Lie algebra $(\fg , \varphi)$
  is a Lie algebra with generators $a^{\beta}_{\alpha}=
  T_{\alpha} a^{\beta}$, $\alpha, \beta \in \ZZ^N$ and
  commutation relations:
  \begin{displaymath}
    [a^{\delta}_{\alpha} , a^{\gamma}_{\beta}]= \delta_{\alpha ,
      \beta \delta^{-1}} a^{\delta \gamma}_{\beta \delta^{-1}} -
    \delta_{\alpha , \beta \gamma} a^{\gamma \delta}_{\beta} \, .
  \end{displaymath}

  Fix $\bar{q} = (q_1 , \ldots, q_N) \in (\CC^\times)^N$ and
  define the homomorphism $\chi_{\bar{q}}: \ZZ^N \to \CC^\times$,
  such that $\chi_{\bar{q}} (\alpha) = \bar{q}^{\alpha}$, where
  $\alpha \in \ZZ^N$ and $\bar{q}^{\alpha}= q_1^{\alpha_1} \cdots
  q^{\alpha_N}_N$.  The Lie algebra $\fs (R, \chi_{\bar{q}})$ of
  formal distributions is a Lie algebra with generators
  $a^{\alpha}_m$, $\alpha \in \ZZ^N$, $m \in \ZZ$.  Define the
  generating functions $a^{\alpha}(z)= \sum_{m \in \ZZ}
  a^{\alpha}_m z^{-m-1}$.  From (\ref{eq:3.4}) we have:
  \begin{displaymath}
    \left[ a^{\alpha}(z) , a^{\beta}(w)\right] =
       \sum_{\gamma \in \ZZ^N}
       (a^{\alpha}_{(\gamma)} a^{\beta})
       (w) \delta (z - \chi (\gamma) w ) \, ,
  \end{displaymath}
  and for modes:
  \begin{eqnarray*}
    \left[ a^{\alpha}_m , a^{\beta}_n \right]
       &=& \sum_{\gamma \in \ZZ^N}
           \chi (\gamma)^m (a^{\alpha}_{(\gamma)} a^{\beta})_{m+n}\\
       &=& \bar{q}^{- \alpha m}(T_{-\alpha} a^{\alpha + \beta})_{m+n}
           - \bar{q}^{\beta m} (a^{\alpha + \beta})_{m+n} \\
       &=& (\bar{q}^{\alpha n}- \bar{q}^{\beta m})
           a^{\alpha + \beta}_{m+n}
  \end{eqnarray*}
  This is a vector generalization of $\sin$-algebra,
  considered in \cite{G-K-L}.
\end{example}

\goodbreak
\section {Non-commutative generalization of $\Gamma$-local formal
  distributions.}
\label{sec:5}

In \S\ref{sec:3} we associated to a $\Gamma$-conformal algebra $R$ a Lie
algebra $\fs (R, \chi)$ of $\chi (\Gamma)$-local formal distributions, by
fixing a homomorphism $\chi: \Gamma \to \CC^\times$.  Conversely, each Lie
algebra of $\chi(\Gamma)$-conformal formal distribution defines a
right $\chi(\Gamma)$-conformal algebra.  In this paragraph we will
consider more general case, where $\chi(\Gamma)$ is a subgroup of
the group of meromorphic transformations of $\CC$.  In all
examples we will consider $\chi: \Gamma \to G L_2 (\CC)$.  For $\alpha
\in \Gamma$ we will denote by $\alpha(z)$ the corresponding under
the homomorphism $\chi$ transformation from $G L_2(\CC)$.  Consider
the generalization of the property (\ref{eq:1.2}) of the
$\delta$-function:

\begin{proposition}
  \label{prop:5.1}
  Let $\alpha(z)$ and $\beta(z)$ be from $G L_2(\CC)$, then:
  \begin{eqnarray}
    \label{eq:5.1}
    \alpha' (z) \delta (\alpha(z)- w)
       &=& \delta(z- \alpha^{-1} (w)) \\
    \nonumber
    \beta'(w) \delta (z- \beta (w))
       &=& \delta (\beta^{-1}(z)- w)
  \end{eqnarray}
\end{proposition}

\begin{proof}
  The group $G L_2(\CC)$ is generated by the transformations a)
  $z \mapsto  \, \alpha z$, $\alpha \in \CC$; b) $z \to z-
  \alpha$, $\alpha \in \CC$; c) $z \mapsto \frac{1}{z}$.  It is
  sufficient to prove (\ref{eq:5.1}) for (a), (b) and (c)
  transformations.  For (a) this is (\ref{eq:1.3}).  To prove
  (b), let
  \begin{mathletters}
    \label{eq:5.2}
    \begin{eqnarray}
      \label{eq:5.2a}
      \delta(z- (w - \alpha)) &=& \sum_{k \in \ZZ} z^{-k-1}
      (w - \alpha)^k \\
      \label{eq:5.2b}
      \delta ((z+ \alpha) - w) &=& \sum_{k \in \ZZ}
      (z+ \alpha)^{- k -1} w^{k} \, .
    \end{eqnarray}
  \end{mathletters}

We define
\begin{equation}
  \label{eq:5.3}
    (z+ w)^{- \ell-1}= \sum^{\infty}_{k=0} \frac{(k+ \ell)!}{k!
      \ell !} z^{- \ell -k-1}(- \alpha)^k
\end{equation}
for $\ell \ge 0$.  Comparing coefficients of $z^{-k-1}
w^{\ell}$ in right-hand sides of (\ref{eq:5.2}) and
(\ref{eq:5.3}) we get (b).  So, the function
\begin{equation}
  \label{eq:5.4}
  \delta (z- w + \alpha) = \sum_k z^{-k-1}(w - \alpha)^k
  = \sum_{\ell} (z+ \alpha)^{- \ell-1} w^{\ell}
\end{equation}
is well defined.  The proof of (c) is similar.
\end{proof}

Further we will suppose for simplicity that $\Gamma \subset G L_2
(\CC)$.  For an $N$-element subset $S \subset \Gamma$ we shall
use the notation
\begin{displaymath}
  (z- w)^N_S = \prod_{\alpha \in S} (z- \alpha (w)) \,.
\end{displaymath}

\begin{definition}
  \label{def:5.1}
  Two formal distributions $a(z)$ and $b(z)$ with coefficients in
  a Lie algebra $\fg$ are called $S$-local if in $\fg[[z,z^{-1} ,
  w , w^{-1}]]$ one has:
\begin{displaymath}
  (z- w)^N_S [a(z), b( w)]=0 \, .
\end{displaymath}
\end{definition}

%%%%%%%%%%%%%%%%%%%%%%%%%%%%%%%%%%%%%%%%%%%%%%%%%%%%%%%%%%%%%%%%%%%%%%%%%%%
%%%%%%%%%%%%%%%%%%%%%%%%%%%%    Jan's part     %%%%%%%%%%%%%%%%%%%%%%%%%%%%
%%%%%%%%%%%%%%%%%%%%%%%%%%%%%%%%%%%%%%%%%%%%%%%%%%%%%%%%%%%%%%%%%%%%%%%%%%%

\begin{proposition}
  \label{prop:5.1a}
  If $a (z)$ and $b (z)$ are $S$-local formal distributions, then
  there exists a unique decomposition
  \begin{equation}
    \label{eq:5.5}
    [a(z), b(w)] = \sum_{\alpha \in S}
      \left( a_{(\alpha^{-1})} b \right) (w) \delta (z - \alpha(w))
  \end{equation}
  The formal distributions $(a_{(\alpha^{-1})} b) (w)$ are given
  by the formula:
  \begin{displaymath}
    \left( a_{(\alpha^{-1})} b \right) (w)
       = \Res_{z} \prod_{\beta \in S \atop \beta \neq \alpha}
         \frac{z - \beta (w)}{\alpha (w) - \beta (w)}
         [a (z), b (w)] \, .
  \end{displaymath}
\end{proposition}

Proof is the same as the proof of Proposition~\ref{prop:1.1}.

For $\alpha \in \Gamma$ introduce the following operator:
\begin{equation}
  \label{eq:5.6}
  T_{\alpha^{- 1}} a (z) = \alpha' (z) a(\alpha (z))
\end{equation}
It is clear, that $T_\alpha$, $\alpha \in \Gamma$, preserve
$\Gamma$-locality.

For $\alpha$-products $(a _{(\alpha)} b) (w)$ defined as \req{5.5} we have
all the properties of left $\Gamma$-conformal algebras (see
Definition~\ref{def:3.1}).  Conversely, if we have a
$\Gamma$-conformal algebra $R$ and a homomorphism $\chi: \Gamma
\to G L_2 (\CC)$, we can construct a Lie algebra $\fs (R, \chi)$
of $\chi (\Gamma)$ local formal distributions as follows.
Consider a vector space over $\CC$ with the basis $a_n$, $a \in
R$, $n \in \ZZ$, and denote by $\fs (R, \chi)$ the quotient of
this space by the $\CC$-space of elements of the form:
\begin{eqnarray}
  &(\lambda a + \beta b)_n - \lambda a_n - \beta b_n \, , &
     \quad
     \lambda, \beta \in \CC \, ,
     \;
     a, b \in R
     \nonumber \\
  &(T_\alpha^{-1} a)_n - (\alpha' (z) a (\alpha (z)))_n, &
     \quad
     \alpha \in \Gamma
  \label{eq:5.7}
\end{eqnarray}
where $(\alpha' (z) a (\alpha (z)))_n$ denotes the coefficient of
$z^{-n-1}$ in the Fourier series decomposition.  Then the formula
\begin{equation}
  \label{eq:5.8}
  [a (z), b (w)] = \sum_{\alpha \in \Gamma} (a_{(\alpha^{-1})} b)
     (w) \delta(z - \alpha (w))
\end{equation}
gives a Lie algebra structure on $\fs (R, \chi)$ of
$\Gamma$-local formal distributions $a (z) = \sum_{n \in \ZZ} a_n
z^{-n-1}$.

\begin{remark}
  \label{rem:5.1}
  If we let $a (z) T_\alpha = \alpha' (z) a (\alpha (z))$, then
  we have all properties of right $\Gamma$-conformal algebras.  We
  prefer (the more customary) left $\Gamma$-modules.
\end{remark}

Consider some important examples.

\begin{example}
\rm
  \label{ex:5.1}
  \emph{Lie algebra of pseudodifferential operators on the
    circle.}

  Let $R$ be the general $\ZZ$-conformal algebra $g c_1 (\ZZ)$
  defined by \req{4.7}.  Fix a homomorphism $\chi: \ZZ \to G L_2
  (\CC)$ defined as $\chi (n) = \left( 1 \quad n \atop 0 \quad 1
  \right)$, $n \in \ZZ$.  The corresponding Lie algebra $\fs
  (\ZZ, \chi)$ is a Lie algebra with basis $a^m_k$, $m, k \in
  \ZZ$.  Let $a^m (z) = \sum_k a^m_k z^{-k-1}$ be the generating
  function.  Then from \req{5.8} we have:
  \begin{eqnarray}
    [a^m (z), a^n (w)]
       &=& \sum_{s \in \ZZ} (a^m_{(s)} a^n) (w)
           \delta (z - w + s)
           \nonumber \\[-.75ex]
    \label{eq:5.9} \\[-.75ex]
       &=& (T_{-m} a^{m+n}) (w) \delta (z - w - m)
           - a^{m+n} (w) \delta (z - w + n)
           \nonumber
  \end{eqnarray}
  To write down explicitly the commutation relations for
  coefficients of formal distributions $a^m (z)$, consider the
  Lie algebra $\pdiff (S^1)$ of pseudodifferential operators on
  the circle (see \cite{Kh-L-R} for more details).  This is a Lie
  algebra of the associative algebra with the basis $x^m
  \partial^n$, $m, n \in \ZZ$.  The commutation relations are
  defined from the commutation relations on $x$ and $\partial$:
  \begin{eqnarray*}
    \partial \circ f(x)      &=& f (x) \circ \partial + f' (x) \\
    \partial^{-1} \circ f(x) &=& \sum^{\infty}_{n = 0}
         (-1)^n f^{(n)} (x) \partial^{-n-1} \, ,
  \end{eqnarray*}
  where $f (x) \in C^\infty (S^1)$.  Introduce a new basis of
  $\pdiff (S^1)$ with generators $x^m D^k$, $m,k \in \ZZ$, where
  $ D = x \partial$.  These generators have more simple
  commutation relations:
  \begin{equation}
    \label{eq:5.10}
    [x^m D^k, x^n D^\ell]
       = x^{m+k} \left((D+n)^k D^\ell - (D + m)^\ell D^k \right)
  \end{equation}
  For $k < 0$ we will understand $(D + n)^k$ as an expansion by
  the negative power of $D$ (see \req{5.3}).  Introduce the
  generating function:
  \begin{equation}
    \label{eq:5.11}
    a^m (z) = - \sum_{k \in \ZZ} x^{-m} D^k z^{-k-1} \, .
  \end{equation}
  Using the property \req{5.4} of $\delta$-function, we have
  \begin{eqnarray}
    \lefteqn{[a^m (z), a^n (w)] } \nonumber\\
       &=& \sum_{k, l \in \ZZ} x^{- (m+n)}
           \Bigl( (D-n)^k D^\ell z^{-k-1} w^{- \ell-1}
%%%                \nonumber \\
%%%       && \quad {}
           - (D - m)^\ell D^k z^{-k-1} w^{-\ell-1}
           \Bigr)
           \nonumber \\
       &=& x^{- (m+n)}
          \biggl( \sum_{k, \ell \in \ZZ} D^k D^\ell (z+n)^{-k-1}
               w^{-\ell -k -1} w^k
%%%               \nonumber \\
%%%       && \quad {}
               - D^\ell D^k (w + m)^{-\ell-k-1}
               (w + m)^k z^{-k-1}
          \biggr)
          \nonumber \\
       &=& a^{m+n} (w + m) \delta (z - w - m)
%%%          \nonumber \\
%%%       && \quad {}
          - a^{m+n} (w) \delta (z - w - n)
          \, .
     \label{eq:5.12}
  \end{eqnarray}
  Comparing \req{5.9} and \req{5.12}, we see that the Lie algebra
  $\fs (\Gamma, \chi)$ is isomorphic to the Lie algebra $\pdiff
  (S^1)$.

  Thus, we have shown, that the Lie algebra of $q-\pdiff (S^1)$
  and the Lie algebra of $\pdiff (S^1)$ correspond to the same
  general $\ZZ$-conformal algebra $g c_1 (\ZZ)$, but to the
  different homomorphisms $\ZZ \to G L_2 (\CC)$.

  More generally, we can consider the general conformal algebra
  $g c_1 (G L_2 (\CC))$ with the usual action of $G L_2
  (\CC)$ on $\CC$.  Let $A, B, C \in G L_2
  (\CC)$.  By \req{4.7} we have:
  \begin{displaymath}
    a^A_{(C)} a^B
      = \delta_{C, A^{-1}} T_{A^{-1}} a^{A B} - \delta_{C, B}
         a^{B A} \, .
  \end{displaymath}
  It is clear, that the general conformal algebra admit reduction to
  the subalgebra.  This reduction is well defined at the level of
  the Lie algebra $\fs (g \ell_1 (G L (\CC)))$:
  \begin{equation}
    \label{eq:5.13}
    [a^A (z), a^B (w)]
       = A (w)'_w a^{A B} (A (w))
          \delta (z - A (w)) - a^{B A} (w) \delta
          (z - B^{-1} (w)) \, .
  \end{equation}
  As was shown, for the subgroup $H_1 = \left\{ \Bigl( {
      q^{\frac{m}{2}} \quad\;\; 0 \;\; \atop \quad 0 \quad\;\;
      q^{- \frac{m}{2}} } \Bigr), \; m \in \ZZ \right\}$ we get
  the Lie algebra of $q - \pdiff (S^1)$; for the subgroup $H_2 =
  \left\{ \left( 1 \quad m \atop 0 \quad 1 \right), \; m \in \ZZ
  \right\}$ we get the Lie algebra $\pdiff (S^1)$.  For the
  subgroup $H_3 \subset G L_2 (\CC)$, generated by two elements
  $a = \left( 1 \quad 1 \atop 0 \quad 1 \right)$ and $b = \Bigl(
  { q^{\frac{1}{2}} \;\quad\; \atop \;\quad\; q^{-\frac{1}{2}} }
  \Bigr)$ we get the Lie algebra, that contains as subalgebra the
  Lie algebra of $\pdiff (S^1)$ and $q - \pdiff (S^1)$.  Taking
  different subgroups $H$ in $G L (2, \CC)$ (not necessarily
  discrete) we get a family of infinite-dimensional Lie algebras
  $\fs (g c_1 (H))$.
\end{example}

\begin{example}
\rm
  \label{ex:5.2}
  Consider the group $\Gamma$ generated by two elements $\epsilon$ and $T$
  with relations: $\epsilon^2 = 1$ and $\epsilon T = T^{-1}
  \epsilon$.  Fix two homomorphisms $\varphi_{i}: \Gamma \to
  \Aut (g \ell_\infty)$, $i = 1, 2$.
  \begin{displaymath}
    \varphi_1 : \left\{
      \begin{array}{l}
        \epsilon (E_{i j}) = - E_{-j, -i} \\
        T (E_{i j}) = E_{i + 1, j + 1}
      \end{array}
    \right.
    \hbox{\quad and \quad}
    \varphi_2: \left\{
      \begin{array}{l}
        \epsilon (E_{i j}) = E_{- i, -j} \\
        T (E_{i j}) = E_{i + 1, j + 1}
      \end{array}
    \right.
  \end{displaymath}
  In both cases $g \ell_\infty$ is a free $\CC [\Gamma]$ module
  with the basis $A^m = E_{0, m}$, $m \in \ZZ_+$.  The
  $\Gamma$-products on generators $A^m$ are as follows:

  \alphalist
  \begin{enumerate}
  \item %% a)
    For $\varphi_1$ we have:
    \begin{eqnarray*}
      A^m_{(\gamma)} A^n &=& [T^r A^m, A^n]
         = \delta_{r, -m} T^{-m} A^{m+n}
            -\delta_{r, n} A^{m+n}
%%            \\
            \quad
%%            && \qquad
            \hbox{for }
            \gamma = T^r \, , \\
      A^m_{(\gamma)} A^n &=& [\epsilon T^r A^m, A^n]
         = -\delta_{r, 0} T^{-m} A^{n+m}
            +\delta_{m+n, -r} A^{m+n}
%%            \\
            \quad
%%            && \qquad
            \hbox{for }
            \gamma = \epsilon T^r \, .
    \end{eqnarray*}

  \item %% b)
    For $\varphi_2$ we have:
    \begin{eqnarray*}
      A^m_{(\gamma)} A^n
        &=& \delta_{r, -m} T^{-m} A^{m+n}
            -\delta_{r, n} A^{m+n}
            \quad
            \hbox{for }
            \gamma = T^r \, , \\
      A^m_{(\gamma)} A^n
        &=& \delta_{r, -m} T^{m} B^{n-m}
            -\delta_{r, -n} B^{n-m}
            \quad
            \hbox{for }
            \gamma = \epsilon T^r \hbox{ and } n \ge m \, .
    \end{eqnarray*}
    These give us two structures of a $\Gamma$-conformal algebra,
    which we will denote $R_1$ and $R_2$ respectively.
    Consider two homomorphisms $\chi_i : \Gamma \to G L_2 (\CC)$,
    $i = 1, 2$:
    \begin{eqnarray*}
      \hbox{a)} \quad
         \chi_1 (\epsilon) &=& \left(
           \begin{array}{cc}
             0   &  1  \\
             -1  &  0
           \end{array}
         \right) \, ,
         \quad
         \epsilon_1 (z) = - \frac{1}{z} \, ; \\
%%%         \quad
         \chi_1 (T) &=& \Bigl(
         \begin{array}{cc}
           q^{\frac{1}{2}}  &  0  \\
           0                &  q^{-\frac{1}{2}}
         \end{array}
         \Bigr) \, ,
         \quad
         T_1 (z) = q z \, ,
         \\
      \hbox{b)} \quad
         \chi_2 (\epsilon) &=&
         \left(
           \begin{array}{cc}
             1   &  0  \\
             0   &  -1
           \end{array}
         \right) \, ,
         \quad
         \epsilon_2 (z) = - {z} \, ; \\
%%%         \quad
         \chi_2 (T) &=&
         \left(
           \begin{array}{cc}
             1  &  1  \\
             0  &  1
           \end{array}
         \right) \, ,
         \quad
         T (z) = z + 1 \, ,
    \end{eqnarray*}
    We can define four Lie algebras $\fs (R_i, \chi_j)$, $i, j =
    1, 2$, of $\Gamma$-local formal distributions with generators
    $A^m_n$, $m \in \ZZ_+$, $n \in \ZZ$, and the generating
    function $A^m (z) = \sum_n A^m_n z^{-n-1}$.

    \alphalist
    \begin{enumerate}
    \item %%a)
      $S (R_, \chi_1)$:
      \begin{eqnarray*}
        \lefteqn{\hspace{-1em}
          [A^m (z), A^n (w)] } \\
           \qquad
           &=& q^m A^{m+n} (q^m w) \delta
               (z - q^m w) - A^{m+n} (w) \delta (z - q^{-n} w)  \\
           && \quad {}
              - q^m A^{m+n} (q^m w) \delta \left( z + \frac{1}{w}
                \right)
              + A^{m+n (w)} \delta \left( z + \frac{q^{m+n}}{w} \right)
      \end{eqnarray*}
      For the coefficients $A^m_k$ we have the commutation
      relations
      \begin{displaymath}
        [A^m_k, A^n_\ell] = (q^{-m \ell} - q^{-n k})
           A^{m+n}_{k+\ell} - (-1)^k q^{m k} (q^{-n \ell}
           - q^{n k}) A^{m+n}_{\ell-k}
      \end{displaymath}
      This is the $C$-series of Lie algebras of quantum
      torus (\cite{G-K-L}).

    \item %%b)
      $\fs (R_2, \chi_1)$:
      \begin{eqnarray*}
        [A^m (z), A^n (z)] &=& A^{m+n} (q^m w) \delta (z - q^n w)
           - A^{m+n} (z - q^{-n} w) \\
         &&\quad{}
           - q^{-m} A^{n-m}
           (q^{-m} w) \delta \left( z + \frac{q^m}{w} \right)
%%           \\
%%         &&\quad{}
           + A^{n-m} (w) \delta \left( z + \frac{q^n}{w} \right) \, ,
           \quad \\
         &&
           \hbox{for }
           m,n \in \ZZ\, ,
           \;
           n \ge m \\
        \noalign{\hbox{\hspace*{50pt} and}}
        [A^m_k, A^n_\ell] &=& (q^{-m \ell} - q^{- k \ell})
           A^{m+n}_{k+\ell} - (-1)^k
           (q^{m \ell} - q^{n k})
           A^{n-m}_{\ell-k} \, .
      \end{eqnarray*}

    \item %%c)
      $S (R_1, \chi_2)$:
      \begin{eqnarray}
        [A^m (z), A^n (w)] = A^{m+n} (w+m) \delta (z - w - m)
%%            \nonumber \\
%%          &&\quad{}
          - A^{m+n} (w) \delta (z - w + n)
            \nonumber \\
          \quad{}
           - A^{m+n} (w + m) \delta (z + w)
%%            \nonumber \\
%%          &&\quad{}
           + A^{m+n} (w) \delta (z + w - m - n)
        \label{eq:5.14}
      \end{eqnarray}
      Consider in $\pdiff (S^1)$ with the basis \req{5.11} a
      subalgebra, stable under the automorphism $w$ defined
      by: $w (x^n D^k) = x^m (m - D)^k$.  This is a Lie
      algebra with the basis $x^m D^k - w (x^m D^k)$.
      Introduce the generating function of the form:
      \begin{displaymath}
        C^m (z) = - \sum_{k \in \ZZ} (x^{-m} D^k - x^{-m} (-m
        -D)^k) z^{-k-1} = a^m (z) + a^m (m - z) \, ,
      \end{displaymath}
      where $a^m (z)$ are given by \req{5.12}.  It is easy to
      check, that the fields $C^m (z)$ satisfy equation
      \req{5.14}.  We will call this subalgebra a $C$-series of
      $\pdiff (S^1)$.

    \item %%d)
      $S (R_2, \chi_2)$ For generating functions we have:
      \begin{eqnarray}
        \label{eq:5.15}
        [A^m (z), A^n (w)] = A^{m+n} (w + m) \delta (z - w - m)
%%            \nonumber \\
%%          &&\quad{}
           - A^{m+n} \delta (z - w + n)
             \nonumber \\
          \quad{}
           + A^{n-m} (w-m) \delta (z + w - m)
%%            \nonumber \\
%%          &&\quad{}
           - A^{n-m} \delta (z + w + n) \
      \end{eqnarray}
      for $m,n \in \ZZ_+$, and $n \ge m$.  This subalgebra is
      called the $B$-series of $\pdiff (S^1)$.  It can be defined
      as a subalgebra stable under the second-order automorphism
      $\sigma$ defined by $\sigma (x^n D^n) = x^{-m} (-D)^n$.
      Relations \req{5.15} are exactly the relations on
      generating series $B^m (z) = \sum_{k \in \ZZ} (x^{-m} D^k +
      \sigma (x^{-m} D^k)) z^{-k-1} = a^m (z) - a^{-m} (-z)$.
    \end{enumerate}

  \end{enumerate}

\end{example}
\subsection*{Acknowledgments}

The research of M.~G.-K. was supported in part
by INTAS grant 942317 and Dutch NWO organization, by RFFI
grant 96-0218046, and grant  96-15-96455  for support of scientific schools.
M.~G.-K. is also very grateful to the Institut Girard
Desarques, URA 746, Universit\'{e} Lyon-1 for the kind hospitality
and support.

\end{document}